\documentclass[a4paper,11pt]{article}

\usepackage[a4paper, portrait]{geometry}
\usepackage[applemac]{inputenc}
\usepackage{amsmath}
\usepackage{slashed}
\usepackage{amscd}\usepackage{amstext}\usepackage{amsbsy}\usepackage{amsopn}\usepackage{amsthm}\usepackage{amsxtra}\usepackage{upref}\usepackage{amsfonts}
\usepackage{amssymb}\usepackage{euscript}
\usepackage[]{latexsym}\usepackage{graphicx}\usepackage{color}\usepackage[all]{xy}
\usepackage[american]{babel}
\usepackage{slashed}
\usepackage{tikz}
\usepackage{tikz-3dplot}
\usepackage{mdframed}
\usepackage{todonotes}
\usetikzlibrary{patterns}
\usepackage{hyperref}
\usepackage{color}
\usepackage{cleveref}

\newtheorem{conjecture}{Conjecture}

\newtheorem{theorem}{Theorem}
\newmdtheoremenv{framedtheorem}{Theorem}

\def\M{\mathbb{M}}

\def\part{\mathcal{P}}

\title{A short introduction to Monstrous Moonshine}
\author{Valdo Tatitscheff \\ \small{valdo.tatitscheff@normalesup.org} \\ \small{IRMA, UMR 7501, Universit\'e de Strasbourg et CNRS}\\
		\small{7 rue Ren\'e Descartes 67000 Strasbourg, France}}
\date{January 24, 2019}

\begin{document}

\maketitle
\abstract{This article is a short and elementary introduction to the monstrous moonshine aiming to be as accessible as possible. I first review the classification of finite simple groups out of which the monster naturally arises, and features of the latter that are needed in order to state the moonshine conjecture of Conway and Norton. Then I motivate modular functions and modular forms from the classification of complex tori, with the definitions of the $J$-invariant and its $q$-expansion as a goal. I eventually provide evidence for the monstrous moonshine correspondence, state the conjecture, and then present the ideas that led to its proof. Lastly I give a brief account of some recent developments and current research directions in the field.}

\section{Finite simple groups and the monster}

\paragraph{Why finite simple groups ?}

A \textit{group} $G$ is a set endowed with an associative inner law with neutral element $e\in G$, and such that each $x\in G$ has an inverse. A \textit{subgroup} $H$ of $G$ is a non-empty subset of $G$ closed under multiplication and inverses, in which case one writes $H\leq G$. Any subgroup $H$ of $G$ defines a partition of $G$ into left and right $H$-orbits; the left (respectively, right) $H$-orbit of a $g\in G$ is written $gH$ (respectively, $Hg$). The set of left (respectively, right) orbits is denoted $G/H$ (resp., $H\backslash G$). Left and right $H$-orbits of a $g\in G$ need not coincide. However, there are as many left $H$-orbits in $G$ as right $H$-orbits. This number is called the \textit{index} of $H$ in $G$. 

A subgroup $H$ of $G$ is said to be \textit{normal} (and one writes $H\unlhd G$) if the following holds:
\begin{equation}
	\forall g\in G,\ gH=Hg \ .
\end{equation}
In that case, the set $G/H\simeq H\backslash G$ can be endowed with a group structure. The \textit{group} $G/H$ is called the \textit{quotient group} of $G$ by $H$. A group $G$ always has two \textit{trivial} normal subgroups: $\{e\}$, and $G$ itself. These yield the quotient groups $G/\{e\}=G$ and $G/G=\{e\}$.

When $G$ has a non-trivial normal subgroup $H$, the knowledge of $H$ and $G/H$, which are two "smaller" groups, is useful in order to understand the structure of $G$. \textit{Simple groups} are the groups which do not have any non-trivial normal subgroup, and thus cannot be ``broken into smaller pieces", in the following sense.

A \textit{composition series} for a group $G$ is a finite sequence
\begin{equation}
	\{e\}=G_0\unlhd G_1 \unlhd \cdots \unlhd G_n=G \ ,
\end{equation}
where $n\geq 1$ is an integer, such that for any $1\leq i\leq n$ one has $G_{i-1}\unlhd G_i$, and where each \textit{factor} $G_i/G_{i-1}$ is simple. Groups do not always admit composition series ($\mathbb{Z}$ does not admit any, for example), however \textit{every finite group does}. 

Furthermore, a group $G$ may admit different composition series. For instance, the cyclic group $C_{12}$ fits into:

\begin{equation}
	\left.
	\begin{array}{ccccccc}
	\{e\} & \unlhd & C_3 & \unlhd & C_6 & \unlhd & C_{12} \\
	\{e\} & \unlhd & C_2 & \unlhd & C_6 & \unlhd & C_{12} \\
	\{e\} & \unlhd & C_2 & \unlhd & C_4 & \unlhd & C_{12}
	\end{array} \ .
	\right.
\end{equation}

Fortunately, the \textit{Jordan-Hölder theorem} says that any two composition series for a finite group $G$ have the same factors, up to permutation. The set of factors for $C_{12}$, for example, is always $\{C_2,C_2,C_3\}$. 

Two non-isomorphic groups may have the same factors. For instance, and for $p$ a prime number, the dihedral group $D_p$ (the group of symmetries of a regular $p$-gon) and the cyclic group $C_{2p}$ both have exactly $\{C_2,C_p\}$ as their set of factors. 

We may make an analogy with the pursuit of building structures with small plastic bricks. Given a set of bricks of different shapes and different colors, one can built many objects by assembling them together differently just as factors do not characterize a group. 

The Jordan-Hölder Theorem corresponds to the fact that even if there are many ways to build the same object, one must use the same set of bricks. 

Since there often is a catalog of the different bricks in the game, one may wonder if it is possible to write a catalog of the finite simple groups. It has actually been done, and it is one of the greatest mathematical achievements of the last century.

\paragraph{The classification of finite simple groups}

Group theory finds its roots in $1832$, as Galois introduced \textit{permutation groups} and \textit{normal subgroups}, and proved that the groups $A_n$ of even permutations of a set of $n$ elements and $\mathrm{PSL}_2(\mathbb{F}_p)$ are simple respectively when $n\geq5$ and $p\geq5$. However, the abstract definition of groups as the one we gave above was only formulated by Cayley more than twenty years later \cite{cayley1854vii}. The \textit{Sylow theorems}, cornerstones of the theory of finite groups, were published in $1872$ in \cite{sylow1872theoremes}. Twenty more years later, Hölder proves that the order of non-abelian finite simple groups is a product of at least four primes, and asks for a general classification of finite simple groups \cite{holder1892einfachen}. 

The first half of the 20th century saw the development of new ideas in finite group theory. In Note M of the second edition of his book \cite{burnside1911theory}, Burnside conjectures that every non-abelian finite simple group has even order, motivated by numerous observations that odd-order and even-order groups behave very differently. The \textit{Feit-Thompson theorem} \cite{feit1963solvability} eventually completed the proof of Burnside's intuition, proving that odd-order groups are solvable. This opened the doors to study and classify finite simple groups since a consequence of the Feit-Thompson theorem is that every non-abelian finite simple group has an involution. One could thus try to classify centralizers of involutions. In $1972$ Gorenstein announced a program to \textit{achieve the classification of finite simple groups} \cite[Appendix]{gorenstein1979classification}. It was eventually completed in $2004$ by Aschbacher and Smith \cite{aschbacher2004classification} (even if in fact, a small oversight of the classification was corrected in $2008$). For a more detailed account of the history of the classification of finite simple groups, see \cite{solomon,aschbacher2004status} and references therein.

\begin{theorem}[Classification of the finite simple group]
Let $G$ be a finite simple group. Then $G$ is either:
\begin{itemize}
\item an element of one of the 18 infinite families.
\item one out of 26 exceptions called sporadic groups.
\end{itemize}
\end{theorem}

Among the infinite families, one finds the sequences $(C_p)_{p\ \mathrm{prime}}$ of cyclic groups of prime order and $(A_n)_{n\geq5}$ of alternating permutation groups. The 16 other families are the so-called \textit{finite simple groups of Lie type}, and are grouped together according to their structure. For example, nine of these form the \textit{Chevalley groups}:
\begin{equation}
\left(A_n(q)\right)_{n\geq 1},\ \left(B_n(q)\right)_{n\geq 2},\ \left(C_n(q)\right)_{n\geq 3},\ \left(D_n(q)\right)_{n\geq 4},\ 
G_2(q),\ F_4(q),\ E_6(q),\ E_7(q),\ E_8(q) \ ,
\end{equation}
where $q=p^\alpha$ is a positive power of a prime number $p$. The are four additional families of \textit{Steinberg groups}, one of \textit{Suzuki groups} and two of \textit{Ree groups} (to which is appended the \textit{Tits group}).

The 26 sporadic cases are the only finite simple groups which do not belong to any of the infinite families mentioned above. Some of them are relatively small (in terms of their \textit{order}, that is, their cardinal as a set). The smallest is the Mathieu group $M_{11}$ which has $7920$ elements. There are also a few mastodons. The next-to-biggest one is called the \textit{baby monster} $\mathbb{B}$, and it has around $4\times 10^{33}$ elements. It is nevertheless orders of magnitude smaller than the biggest one which is called the \textit{monster} $\M$, and whose order is:

\begin{equation}
	|\M|=2^{46} \cdot 3^{20} \cdot 5^9 \cdot 7^6 \cdot 11^2 \cdot 13^3 \cdot 17 \cdot 19 \cdot 23 \cdot 29 \cdot 31 \cdot 41 \cdot 47 \cdot 59 \cdot 71\sim8\times10^{53}\ .
\end{equation}

\paragraph{The rise of the monster} While studying $3$-transposition groups (which led for example to the discovery of the sporadic groups $\mathrm{Fi}_{22}$, $\mathrm{Fi}_{23}$ and $\mathrm{Fi}_{24}'$) \cite{fischer}, Fischer provided evidence suggesting the existence of the simple finite group now known as the \textit{baby monster} $\mathbb{B}$. Even if the existence of the latter was only proved in 1977 by Leon and Sims \cite{leonsims}, Fischer (in 1973, unpublished) and Griess \cite{griess1976structure} predicted the \textit{monster group} $\M$ as a simple finite group containing a double cover of $\mathbb{B}$. Just as its baby, the monster was only proved to exist about a decade after it had been predicted. However, its conjectural properties were enough to compute its order, as well as others of its features.

Let $G$ be a finite group, and let $\gamma\in G$. The \textit{conjugacy class} of $\gamma$ in $G$ is the set
\begin{equation}
	\mathcal{C}(\gamma)=\left\{g\gamma g^{-1}\vert g\in G\right\}\ .
\end{equation}
The conjugacy classes partition the group $G$. The monster, if any, was expected to have $194$ conjugacy classes. 

A (complex linear finite-dimensional) \textit{representation} of $G$ is a pair $(\rho,V)$ where $V$ is a finite-dimensional complex vector space, and $\rho:G\rightarrow \mathrm{GL}(V)$ is a group morphism. The group $G$ is hence seen as a finite subgroup of the group of invertible $n\times n$ matrices acting on $V$. The dimension of $V$ is by definition the \textit{dimension of the representation}. Every group has a \textit{trivial representation} (where $V$ is a one-dimensional complex vector space and $\rho$ is constant and equal to $1$), and a \textit{regular representation}, whose vector space is of dimension the order of the group, with a basis $(e_g)_{g\in G}$, and where $\rho$ is the morphism assigning to each $h\in G$ the matrix $\rho(h)$ defined on the basis as:

\begin{equation}
	\rho(h)e_g=e_{hg},\ \forall g\in G\ .
\end{equation}

Every representation of a finite group $G$ is isomorphic to a direct sum of special representations called \textit{irreducibles}. What that means, is that up to a change of basis in $V$, one can assume $\rho$ to be valued in block-diagonal matrices whose sizes are the dimensions of irreducible representations. There are as many irreducible representations of a group $G$ as conjugacy classes in $G$. Needless to say, any piece of information concerning those irreducible representations, such as their dimension, is a big step forward in the understanding of $G$. See \cite{serre} or \cite{fultonharris} for a general introduction to representations. 

Going back to the monster, under the assumption that the smallest non-trivial irreducible representation is $196883$-dimensional (which was conjectured from group-theoretic considerations of Conway and Norton), the $194$ irreducible dimensions were conjectured in $1978$ by Fischer, Livingstone and Thorne through clever computer calculations. This $194$-uplet is encoded in the Online Encyclopaedia of Integer Sequences (OEIS) as the sequence $A001379$. The first few terms are:

\begin{equation}
	(r_n)_{n=1..194}=(1,196883,21296876,842609326,18538750076,...) \ .
\end{equation}

\section{Complex tori, rank-2 lattices and modularity}

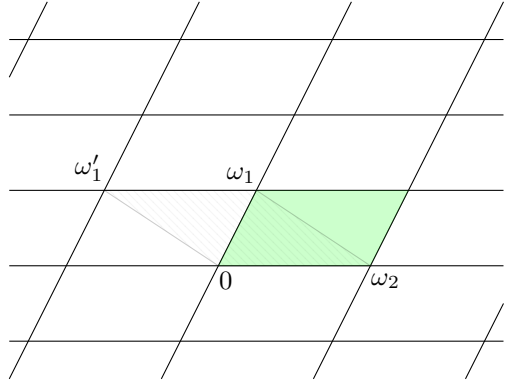
\begin{figure}[h!]
\centering
\begin{tikzpicture}[scale=0.9]
\draw (-0.25,0)--(6.25,0);
\draw (-0.25,1)--(6.25,1);
\draw (-0.25,2)--(6.25,2);
\draw (-0.25,3)--(6.25,3);
\draw (-0.25,4)--(6.25,4);
\draw (-0.25,-0.5)--(2.25,4.5);
\draw (1.75,-0.5)--(4.25,4.5);
\draw (3.75,-0.5)--(6.25,4.5);
\draw (-0.25,3.5)--(0.25,4.5);
\draw (5.75,-0.5)--(6.25,0.5);
\draw (2.6,0.8) node{$0$};
\draw (2.8,2.2) node{$\omega_1$};
\draw (4.7,0.8) node{$\omega_2$};
\draw (0.8,2.3) node{$\omega_1'$};
\draw [fill=green, opacity=0.2] (2.5,1) -- (4.5,1) -- (5,2) -- (3,2) --cycle;
\draw [pattern=north west lines, opacity=0.2] (2.5,1) -- (4.5,1) -- (3,2) -- (1,2) --cycle;
\end{tikzpicture}
\caption{Two bases for the lattice $\Lambda_\tau$.}\label{Fig:lattice}
\end{figure}

\paragraph{Isomorphism classes of complex tori} A \textit{lattice} $\Lambda$ in $\mathbb{C}$ is a subgroup of $\mathbb{C}$ isomorphic to $\mathbb{Z}^2$ such that the real vector space it generates is 2-dimensional (hence it is the whole complex plane). Equivalently, there exists $\omega_1, \omega_2\in\mathbb{C}^\times$ satisfying $\omega_1/\omega_2\notin\mathbb{R}$, and such that
\begin{equation}
	\Lambda=\mathbb{Z}\cdot\omega_1+\mathbb{Z}\cdot\omega_2=\{a\omega_1+b\omega_2|a,b\in\mathbb{Z}\}\subset\mathbb{C} \ .
\end{equation}
Let $\Lambda(\omega_1,\omega_2)$ be the lattice generated by $\omega_1$ and $\omega_2$, shown in \Cref{Fig:lattice}.

A (compact one-dimensional) \textit{complex torus} is a \textit{quotient} $\mathbb{C}/\Lambda$. Topologically $\mathbb{C}/\Lambda(\omega_1,\omega_2)$ is indeed the torus one obtains by gluing together the opposite sides of the parallelogram with vertices $0,\omega_1$ and $\omega_2$, with the same orientation. Gluing the first pair of edges yields a cylinder, and the second, a \textit{torus}. 

There is only one \textit{topological torus}, meaning that any torus can be deformed into any other continuously, that is, without tearing the surface. These deformations are called \textit{homeomorphisms}. \textit{Complex isomorphism} (also called \textit{conformal isomorphisms}) are more constraining than homeomorphisms, since on top of being homeomorphisms, they must also \textit{preserve the angles} at each point. Subsequently, there might be more than one complex isomorphism class of tori up to complex isomorphisms.

One can show that two complex tori $\mathbb{C}/\Lambda_1$ and $\mathbb{C}/\Lambda_2$ are conformally equivalent if and only if there exists an $\alpha\in\mathbb{C}^\times$ such that $\Lambda_2=\alpha\Lambda_1$. Hence one can assume (possibly after multiplying by $\omega_2^{-1}$) that the lattice we consider is of the form $\Lambda_\tau=\Lambda(\tau,1)$ where $\tau=\omega_1/\omega_2$. Moreover (possibly after exchanging $\tau$ and $1$ and then multiplying by $\tau^{-1}$), one can assume that $\tau$ lies in the complex upper half-plane, also called the \textit{Poincaré plane} and denoted $\mathbb{H}$. 

Hence we have defined a complex torus as an equivalence class of lattices in the complex plane, the equivalence relation being the multiplication by non-zero complex numbers. In each equivalence class there is a lattice of the form $\Lambda_\tau$ for some $\tau\in\mathbb{H}$. 

It can still be that  $\Lambda_\tau=\Lambda_{\tau'}$ for two distinct $\tau,\tau'\in\mathbb{H}$ as in \Cref{Fig:lattice}, where the pairs $(\omega_1,\omega_2)$ and $(\omega_1',\omega_2)$ both generate the same lattice. To tackle this issue, let us consider a lattice $\Lambda = \Lambda(\omega_1,\omega_2) = \Lambda(\Omega_1,\Omega_2)$. Then there are $A_i,B_i,a_i,b_i\in\mathbb{Z}$, for $i\in\{1,2\}$, such that:

\begin{equation}
\left\{
\begin{array}{ccc}
\Omega_1 & = & A_1\omega_1+B_1\omega_2 \\
\Omega_2 & = & A_2\omega_1+B_2\omega_2 \\
\omega_1 & = & a_1\Omega_1+b_1\Omega_2 \\
\omega_2 & = & a_2\Omega_1+b_2\Omega_2 \\
\end{array}
\right. \ ,
\end{equation}
which leads to:

\begin{equation}
\left\{
\begin{array}{ccccc}
A_1a_1+B_1a_2 & = & A_2b_1+B_2b_2 & = & 1 \\
A_1b_1+B_1b_2 & = & A_2a_1+B_2b_1 & = & 0\\
A_1a_1+B_1a_2 & = & A_2b_1+B_2b_2 & = & 1 \\
A_1b_1+B_1b_2 & = & A_2a_1+B_2b_1 & = & 0 \\
\end{array}
\right. \ .
\end{equation}

We see that the quadruples $(A_1,B_1,A_2,B_2)$ and $(a_1,b_1,a_2,b_2)$ form matrices with integer coefficients which are mutually inverse (hence their determinant is $\pm1$). In other words, changes of bases in a lattice form the group $\mathrm{SL}_2(\mathbb{Z})$. 

Any basis of a lattice $\Lambda_\tau$ with $\tau\in\mathbb{H}$ is of the form $(a\tau+b,c\tau+d)$ for some quadruple $(a,b,c,d)\in\mathbb{Z}^4$ such that $ad-bc=1$. Now:

\begin{equation}
\mathbb{C}/\Lambda_\tau\simeq\mathbb{C}/\Lambda(\tau,1)\simeq\mathbb{C}/\Lambda(a\tau+b,c\tau+d)\simeq\mathbb{C}/\Lambda_\frac{a\tau+b}{c\tau+d} \ .
\end{equation}

Hence $\tau$ and $\tau'\in\mathbb{H}$ are such that $\Lambda_\tau = \Lambda_{\tau'}$ if and only if there exists $M=\left(\hspace{-0.15cm}\begin{array}{cc} a & b \\ c & d \end{array}\hspace{-0.15cm}\right)\in\mathrm{SL}_2(\mathbb{Z})$ for which
\begin{equation}
	M\cdot \tau=\frac{a\tau+b}{c\tau+d}=\tau' \ ,
\end{equation}
and complex tori modulo complex isomorphisms are in bijection with the orbits of the action of $\mathrm{SL}_2(\mathbb{Z})$ on $\mathbb{H}$. 

One can observe that $-\mathrm{Id}\in\mathrm{SL}_2(\mathbb{Z})$ acts trivially on $\mathbb{H}$, hence one should rather consider the \textit{modular group} $\mathrm{PSL}_2(\mathbb{Z})$ obtained by identifying each matrix with its opposite in $\mathrm{SL}_2(\mathbb{Z})$, instead of $\mathrm{SL}_2(\mathbb{Z})$.

\paragraph{The moduli space $\mathcal{M}_{1}$}
We have seen that the natural action of the modular group on $\mathbb{H}$ is:
\begin{equation}
	\left(\hspace{-0.15cm}\begin{array}{cc} a & b \\ c & d \end{array}\hspace{-0.15cm}\right)\cdot\tau=\frac{a\tau+b}{c\tau+d}.
\end{equation}

\begin{figure}[ht]
\begin{minipage}[b]{0.55\linewidth}
\centering
\begin{tikzpicture}[scale=1.2]
\draw (-2.5,0)--(2.5,0);
\draw (0,-0.3)--(0,4);
\draw (2,0) arc (0:180:2);
\draw[green] (1,1.73) arc(60:120:2);
\draw[color=red] (1,1.73)--(1,4);
\draw[color=red] (-1,1.73)--(-1,4);
\draw[color=red] (1,3)--(0.95,2.95);
\draw[color=red] (1,3)--(1.05,2.95);
\draw[color=red] (1,3.05)--(0.95,3);
\draw[color=red] (1,3.05)--(1.05,3);
\draw[color=red] (-1,3)--(-0.95,2.95);
\draw[color=red] (-1,3)--(-1.05,2.95);
\draw[color=red] (-1,3.05)--(-0.95,3);
\draw[color=red] (-1,3.05)--(-1.05,3);
\draw[fill=black] (0,0) circle(0.05);
\draw[fill=black] (0,2) circle(0.05);
\draw[fill=black] (1,1.73) circle(0.05);
\draw[fill=black] (-1,1.73) circle(0.05);
\draw[black!60!green] (-0.6,1.94)--(-0.5,1.94);
\draw[black!60!green] (0.6,1.94)--(0.5,1.94);
\draw[black!60!green] (-0.58,1.88)--(-0.5,1.94);
\draw[black!60!green] (0.58,1.88)--(0.5,1.94);
\draw (-0.2,0.2) node{$0$};
\draw (0.2,1.8) node{$i$};
\draw (1.3,1.73) node{$\rho$};
\draw (-1.3,1.73) node{$\rho^2$};
\draw (-0.5,3.5) node{$\Delta$};
\draw (2.2,3.5) node{$\mathbb{H}$};
\draw[pattern=north west lines, opacity=0.4] (1,4) -- (1,1.73) -- (1,1.73) arc(60:120:2) -- (-1,1.73) -- (-1,4);
\end{tikzpicture}
$$(\tau\in\Delta)\ \Leftrightarrow\ |\tau|\geq1\ \mathrm{and}\ -\frac{1}{2}\leq\mathrm{Re}(\tau)\leq\frac{1}{2}$$
\end{minipage}
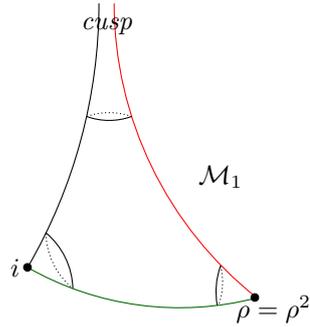
\begin{minipage}[b]{0.35\linewidth}
\centering
\begin{tikzpicture}[scale=1.2]
\draw (-0.1,5) arc (0:-30:7);
\draw[red] (0.1,5) arc (180:231.4:5);
\draw[black!60!green] (-1.04,1.5) arc(240:284.8:4);
\draw[fill=black] (-1.04,1.5) circle(0.05);
\draw[fill=black] (1.95,1.1) circle(0.05);
\draw (-1.2,1.5) node{$i$};
\draw (2.2,0.95) node{$\rho=\rho^2$};
\draw (-0.8,1.96) arc (50:1.7:1);
\draw[densely dotted] (-0.8,1.95) arc (182:230:1);
\draw (1.5,1.52) arc (160:190.4:1);
\draw[densely dotted] (1.5,1.52) arc (14:-20:0.9);
\draw (-0.26,3.5) arc (250:290:0.86);
\draw[densely dotted] (-0.26,3.5) arc (110:70:0.86);
\draw (1.5,2.7) node{$\mathcal{M}_1$};
\draw (0,4.7) node{$\mathrm{\textit{cusp}}$};
\end{tikzpicture}
\label{fig:figure2}
\end{minipage}
\caption{A fundamental domain in $\mathbb{H}$ and the resulting quotient space.}\label{Fig:fundom}
\end{figure}

This action behaves nicely, and in particular it is possible to choose a \textit{fundamental domain} in $\mathbb{H}$, that is, a closed subset of $\mathbb{H}$ which is connected and simply-connected (i.e it has only one component, and no hole), and which contains exactly one point of each orbit under the action of $\mathrm{PSL}_2(\mathbb{Z})$ (except maybe on its boundary). The iterates of this domain under $\mathrm{PSL}_2(\mathbb{Z})$ \textit{tile} the upper-half plane $\mathbb{H}$. The domain $\Delta\subset\mathbb{H}$ shown in \Cref{Fig:fundom} is a usual choice.

The set of orbits $\mathcal{M}_1=\mathbb{H}/\mathrm{PSL}_2(\mathbb{Z})$ has a natural topology induced by the one of $\mathbb{H}$, that one can visualize by gluing the two vertical sides of $\Delta$ together (both going up), as well as the two arcs $(\rho^2,i)$ and $(\rho,i)$ (in that order, meaning that $\rho=e^{2\pi i/3}$ is identified with $\rho^2$). The resulting space has a \textit{cusp} going to infinity, as well as two corners, but if one adds a point at the cusp so the surface closes there, and smoothens the latter at the corners, the resulting space is smooth and homeomorphic to a sphere. This space $\mathcal{M}_{1}$ is called the \textit{moduli space of complex tori} as it is in bijection with isomorphism classes of complex tori. Actually $\mathcal{M}_1$ even inherits a \textit{complex structure} from $\mathbb{H}$.

\paragraph{Modular functions and modular forms}

One may wonder if there are any complex functions on $\mathbb{H}$ invariant under the action of $\mathrm{PSL}_2(\mathbb{Z})$, that is, functions $f$ such that:
\begin{equation}
	f(\tau)=f\left(M\cdot\tau\right)=f\left(\frac{a\tau+b}{c\tau+d}\right)
\end{equation}
for any $M=\mathrm{Mat}(a,b,c,d)\in\mathrm{PSL}_2(\mathbb{Z})$. Since $\mathrm{PSL}_2(\mathbb{Z})$ is \textit{generated} by the two matrices: 
\begin{equation}
	T=\left(\hspace{-0.15cm}\begin{array}{cc} 1 & 1 \\ 1 & 0 \end{array}\hspace{-0.15cm}\right)\ \mathrm{and}\ S=\left(\hspace{-0.15cm}\begin{array}{cc} 0 & -1 \\ 1 & 0 \end{array}\hspace{-0.15cm}\right) \ ,
\end{equation}
meaning that any element in $\mathrm{PSL}_2(\mathbb{Z})$ can be written as a product of $S$'s and $T$'s, the condition for $f$ to be invariant under $\mathrm{PSL}_2(\mathbb{Z})$ becomes:
\begin{equation}
	f(\tau+1)=f(\tau)\ \mathrm{and}\ f\left(\frac{-1}{\tau}\right)=f(\tau) \ .
\end{equation}
Such a (meromorphic) function (in $\mathbb{H}$ and at the cusp) is called a \textit{modular function}. 

In practice, it is difficult to study modular functions directly, and relaxing some of the hypotheses in the definition of modular functions has proved very fruitful. It leads to the definition of modular forms. \\

A \textit{modular form} $f$ of weight $2k$ (for $k\in\mathbb{N}$) is a holomorphic function (in $\mathbb{H}$ and at the cusp) such that:
\begin{equation}
	f(\tau)=(c\tau+d)^{2k}f\left(\frac{a\tau+b}{c\tau+d}\right)
\end{equation}
for any $M=\mathrm{Mat}(a,b,c,d)\in\mathrm{PSL}_2(\mathbb{Z})$. The set $M_k$ of modular forms of weight $2k$ is a complex vector space. Moreover, the product of two modular forms of respective weights $2k$ and $2l$ is a modular form of weight $2(k+l)$, thus
$$M=\bigoplus_{k\geq0}M_k$$
is a ring with unit (which is the constant function equal to $1$).

Non-trivial examples of modular forms can be constructed easily. For example, let $k\in\mathbb{N}_{\geq2}$. Then the $k$-th Eisenstein series defined as:
\begin{equation}
	G_{2k}(\tau)=\sum_{\omega\in\Lambda_\tau\backslash\{0\}} \frac{1}{\omega^{2k}} \ .
\end{equation}
is an example of a modular form of weight $2k$.

The theory of modular forms proves that $M$ is the polynomial ring in $G_4$ and $G_6$:
\begin{equation}\label{ringmod}
	M = \mathbb{C}\left[G_4,G_6\right] \ .
\end{equation}

Since modular functions are not required to be holomorphic but only meromorphic, one can look for modular functions as quotients of modular forms of the same weight, on top of the elements of $M_0$ which are of course modular functions, but only the constant ones by Equation \ref{ringmod}. Since $12$ is the smallest $k$ for which the dimension of $M_{k}$ is strictly bigger than $1$, the simplest non-trivial modular functions can be build out of elements of $M_{12}$. Klein's $J$-invariant is defined as the following ratio of linearly independent elements of $M_{12}$:
\begin{equation}
	J(\tau)=1728\frac{g_2(\tau)^3}{g_2(\tau)^3-27g_3(\tau)^2} \ ,
\end{equation}
where $g_2=60G_4$ and $g_3=140G_6$. It is hence is a modular function, and an important one since \textit{any modular function is a rational fraction of $J$}. Equivalently one can say that the ring (field) of modular functions is $\mathbb{C}(J)$.

\paragraph{$q$-expansions}

Both modular forms and functions must satisfy $f(\tau+1)=f(\tau)$, and hence admit a Fourier decomposition in terms of $q=e^{2\pi i\tau}$. Under $\tau\mapsto q(\tau)$, the upper half-plane $\mathbb{H}$ is mapped to the punctured unit disk, and of course any function $f$ defined as a series of $q$ is invariant under $\tau\mapsto \tau+1$.

It turns out that the coefficients of the Fourier series of some important modular objects (functions or forms) are integers. For example, let $k\in\mathbb{N}_{\geq2}$. Then:
\begin{equation}
	G_{2k}(\tau)=2\zeta(2k)+2\frac{(2\pi i)^{2k}}{(2k-1)!}\sum_{n=1}^\infty \sigma_{2k-1}(n)q^n \ ,
\end{equation}
where $\sigma_k(n)=\sum_{d|n}d^k$. The coefficients of the $q$-expansion of $J$ are also integers, but it is not obvious to understand what they count. These coefficients, which are interesting already because of their definition, are listed in the sequence in referred to as A000521 in the OEIS. The first few terms are:
\begin{equation}
J(\tau)=\sum_{n=-1}^\infty c(n)q^n=\frac{1}{q}+744+196884q+21493760q^2+864299970q^3+20245856256q^4+...
\end{equation}

In the sequel, we will denote $\tilde J$ for the \textit{normalised J-function}: $\tilde J = J-744$. More details concerning the theory of modular functions and modular forms can be found in \cite{milne}, and in the references therein.

\section{Monstrous Moonshine}

\paragraph{First observations and Thompson's conjecture}

Even if the $J$-invariant and the monster group do not seem to be related at all at first glance, they actually are. This surprising correspondence has been dubbed the \textit{Monstrous Moonshine}. In began in 1978 with McKay's famous observation that
\begin{equation}\label{Eq:McKay}
	c(2)=r_1+r_2 \ ,
\end{equation}
where we keep the notations of the previous sections. Explicitly, Eq. \ref{Eq:McKay} is:
\begin{equation}
	196884=196883+1 \ ,
\end{equation}
and it is the huge numbers involved which makes this equality thrilling, as it seems more difficult for it to be a mere coincidence. This first observation was then supplemented by Thompson in \cite{thompson}, who pointed out that:
\begin{equation}
	\left.
	\begin{array}{cc rcrcrcrcrcrcrc}
	c(3) & = & r_1 & + & r_2 & + & r_3 & & & & & & & & \\
	c(4) & = & 2r_1 & + & 2r_2 & + & r_3 & + & r_4 & & & & & & \\
	c(5) & = & 3r_1 & + & 3r_2 & + & r_3 & + & 2r_4 & + & r_5 & & & & \\
	c(6) & = & 4r_1 & + & 5r_2 & + & 3r_3 & + & 2r_4 & + & r_5 & + & r_6 & + & r_7 \\
	\end{array}
	\right. \ ,
\end{equation}
and subsequently made the following conjecture.

\begin{conjecture}[Thompson's conjecture]
There exists a somehow natural $(\natural)$ infinite-dimensional graded representation $\left(\rho_\natural,V^\natural={\bigoplus_{i\geq-1}V^\natural_i}\right)$  of the monster group such that each graded part $V^\natural_i$ is finite dimensional, and such that the generating series of the dimension of these graded parts is the $q$-expansion of the normalized $J$-invariant:
\begin{equation}
	\tilde J(\tau)=\sum_{i\geq-1}\dim \left(V^\natural_i\right) q^i \ .
\end{equation}
\end{conjecture}

\paragraph{Conway-Norton's Monstrous Moonshine}

Let $(\rho, V)$ be a representation of a finite group $G$. The \textit{character} of $(\rho,V)$ is the complex valued function on $G$ defined as:
\begin{equation}
	\left.\begin{array}{rcl}\chi_{(\rho,V)}:G & \rightarrow & \mathbb{C} \\
	g & \mapsto & \mathrm{Tr}\left(\rho(g)\right)\end{array}\right..
\end{equation}

Since $\rho$ is a group morphism and since the trace is cyclically invariant, the character $\chi_{(\rho,V)}$ is constant on each conjugacy class of $G$. Such complex functions on $G$ which are constant on the conjugacy classes are called \textit{class functions}, hence we just saw that any character is a class function.

Let $\mathrm{CC}(G)$ be the set of conjugacy classes in $G$. The vector space of class functions on $G$ is thus $\mathbb{C}^{\mathrm{CC}(G)}$, and it is endowed with a natural hermitian product given by averaging over $G$. The following result is fundamental to the theory of characters of finite groups.

\begin{theorem}
\begin{enumerate}
\item The characteristic functions of the conjugacy classes of $G$ form an orthogonal basis of $\mathbb{C}^{\mathrm{CC}(G)}$.
\item The characters of the irreducible representations of $G$ form an orthonormal basis of $\mathbb{C}^{\mathrm{CC}(G)}$. 
\end{enumerate} 
\end{theorem}

The transition matrix from the basis of characteristic functions to the basis of irreducible characters is called the \textit{character table} of $G$. Computing the character table of a given finite group $G$ is the holy grail of the character theory of $G$. Standard references on finite group theory and characters of finite groups are \cite{serre} and \cite{fultonharris}.\\

Since the dimension of a representation is the trace of $\rho(e)$, Thompson also suggested to study the series
\begin{equation}\label{Eq:McKayThompson}
	T_{[g]}=\sum_{i\geq-1}\mathrm{Tr}\left(\rho_\natural(g)_{\left\vert V^\natural_i\right.}\right) q^i=\frac{1}{q}+\sum_{n=0}^\infty H_n([g])q^n 
\end{equation}
which are now known as \textit{McKay-Thompson series}; there is one for each conjugacy class $[g]$ of $\M$. In Eq. \ref{Eq:McKayThompson}, $g\in\M$ is any representative of $[g]$. The class-functions $H_n$ are called the \textit{Head characters} of $\M$.

We already mentioned the calculation Fischer, Livingstone and Thorne performed in 1978. In fact, they did not only conjectured the dimensions of the irreducible representations of $\M$ but the whole character table of the latter. Building on this Conway and Norton then extended Thompson's conjecture in \cite{conwaynorton} to the following. 

\begin{conjecture}[Conway-Norton's Monstrous Moonshine]
There exists a somehow natural infinite-dimensional graded representation $\left(\rho_\M,V^\natural=\bigoplus V^\natural_i\right)$ of the monster group, such that each graded part $V^\natural_i$ is finite dimensional, and such that for each conjugacy class $[g]$ the McKay-Thompson series $T_{[g]}$ is the $q$-expansion of the normalized Hauptmodul of a subgroup $\Gamma_{[g]}$ of $\mathrm{PSL}_2(\mathbb{R})$ commensurable with $\mathrm{PSL}_2(\mathbb{Z})$.
\end{conjecture}

This requires some explanations. Two subgroups $G$ and $H$ of $\mathrm{PSL}_2(\mathbb{R})$ are said to be \textit{commensurable} if their intersection has finite index in both $G$ and $H$. As an example of a group $\Gamma_{[g]}$, when $[g]=\{e\}$ one has $\Gamma_{\{e\}}=\mathrm{PSL}_2(\mathbb{Z})$, the functions on $\mathbb{H}$ which are $\mathrm{PSL}_2(\mathbb{Z})$-invariant are the modular functions, and we have seen that they can all be expressed as rational functions of $J$ (or, equivalently, of $\tilde J$). This property of $\tilde J$, together with the fact that it is normalised, makes it the $\textit{normalised Hauptmodul}$ for $\mathrm{PSL}_2(\mathbb{Z})$. Some subgroups of $\mathrm{PSL}_2(\mathbb{R})$ commensurable with $\mathrm{PSL}_2(\mathbb{Z})$ also have normalised Hauptmoduln defined in the same way: each modular function for that group is a rational function of the Hauptmodul. The groups for which this is the case are said to be of $\textit{genus zero}$.\\ 

Conway and Norton not only made this conjecture of theirs, but also computed a convincing amount of evidence for it. Using Thompson's work as well as the character table of $\M$ they were able to determine empirically all the groups that should appear in the correspondence (if valid). For example, the data produced by Thompson was enough for them to guess that the class of $\M$ denoted $2A$ would have to correspond to the \textit{Hecke congruence subgroup of level $2$} of $\mathrm{PSL}_2(\mathbb{Z})$. In the notation $2A$, $2$ means that it is a conjugacy class of elements of order $2$ and $A$ distinguishes such classes; there is for instance another class of involutions judiciously called $2B$.

Rather than the $194$ expected McKay-Thompson series one naively expects, there are merely $172$ different ones since for any $g\in\M$, the classes $[g]$ and $[g^{-1}]$ define the same series. Moreover, there are linear relations between the series which even lowers down to $163$ the dimension of the space that the McKay-Thompson series generate. 

Conway and Norton also computed $\mathrm{Tr}\left(\rho_\M(g)_{\left \vert V^\natural_i \right.}\right)$ for all $g\in\M$ and $i\leq10$, and checked that the functions 
\begin{equation}
	\left.\begin{array}{ccc} G & \rightarrow & \mathbb{C} \\
	g & \mapsto & \mathrm{Tr}\left(\rho_\M(g)_{\left\vert V^\natural_i\right.}\right)\end{array}\right.
\end{equation}
are indeed characters of representations of $\M$.

It is in their article \cite{conwaynorton} that the term \textit{monstrous moonshine} was first used to name this outlandish proposed correspondence between the $J$-invariant and the monster group $\M$.

\paragraph{The monster vertex operators algebra}

Even if the monster has been conjectured in the early 1970's, its uniqueness was only proved in $1979$ by Thompson \cite{thompson2} under some hypotheses (a complete proof was only provided in 1989 by Griess, Meierfrankenfeld and Segev \cite{griess1989uniqueness}) and its existence in $1982$ by Griess \cite{griess}, who constructed it explicitly as the automorphism group of a $196884$-dimensional commutative non-associative algebra now called the \textit{Griess algebra}. 

In $1979$, Atkin, Fong and Smith performed a complicated computer calculation \cite{smith} to show that the Head characters $H_n$ defined in Eq. \ref{Eq:McKayThompson} were really characters of the monster. Using powerful results of finite group theory such as \textit{Brauer's characterization of characters}, they could reduce the infinite number of checks this conjecture demands, to a finite one. 

At that time, Frenkel, Lepowsky and Meurman were developing a different approach to $\M$ which uses \textit{vertex operators representations of affine Lie algebras} \cite{frenkellm2}. It would take us too far to enter this construction in details, and instead we refer to \cite{frenkellm} for more on vertex operator algebras. This construction strongly relies on the existence of very symmetrical \textit{even unimodular lattices} $L_{E_8}\simeq\mathbb{Z}^8$ and $\Lambda\simeq\mathbb{Z}^{24}$, respectively the \textit{root lattice of $E_8$} and the \textit{Leech lattice} (the latter was already a cornerstone of Griess' construction of the monster). Both lattices are related to the $q$-expansion of $J$, and furthermore they lead to a particular \textit{vertex operator algebra} whose automorphism group is exactly $\M$. This \textit{monster vertex algebra} is graded and the graded-dimension is exactly the $q$-expansions of $\tilde{J}$. It contains a version of the Griess algebra as its second graded part. Hence the monster module whose existence had been confirmed by Atkin, Fong and Smith was then constructed explicitly.

\paragraph{Borcherd's proof of the monstrous moonshine}
 
Vertex algebras play a very central role in a specific field of theoretical physics called \textit{string theory}. This theory first appeared under the name \textit{dual resonance theory} as an attempt to explain strong interactions, which are responsible for the binding of the nucleons inside the nuclei of atoms, before being demoted to the benefit of another theory called \textit{quantum chromodynamics}. However, around 1974 it was realized that string theory could be an interesting approach to quantum gravity. 

In string theory, what is usually called a particle is a vibration mode of a tiny string which propagates through space-time. The geometry of the latter influences the vibration modes of the string, and \textit{vertex operators} as well as \textit{chiral algebras} (the physics counterpart of vertex operator algebras) encode this dependence. The monster vertex algebra of Frenkel, Lepowsky and Meurman can be interpreted in physics terms as the vertex operator algebra of a \textit{bosonic string theory} in a specific toroidal space-time built out of the Leech lattice. The introduction of \cite{frenkellm} narrates the early years of string theory, and explains the links between string theory and the monstrous moonshine. For a glance into the physical ideas behind string theory a standard reference is the first chapter of the textbook \cite{greschwawitt}. 

Using the \textit{no-ghost theorem} of string theory, Borcherds defined and studied in \cite{proborch} the space of physical states of the string theory associated with the monster vertex operator algebra. He dubbed it the \textit{monster Lie algebra}, and showed that it has the structure of a generalized Kac-Moody algebra. As such, the monster Lie algebra satisfies \textit{twisted denominator identities} which provides enough information on the McKay-Thompson series to fully determine them. This allowed Borcherds to conclude the proof of the monstrous moonshine conjecture of Conway and Norton, which contributed to his Fields medal award in 1998. 

Borcherds' proof and the tools he developed have shed much light both on the structure of $\M$ and on modular forms and functions. For instance, the \textit{Koike-Norton-Zagier formula}
\begin{equation}\label{Eq:KNZ}
p^{-1}\prod_{m,n\in\mathbb{Z},\ m>0} \left(1-p^mq^n\right)^{c(mn)}=J(\sigma)-J(\tau)  
\end{equation}
discovered in the 1980's appears naturally in Borcherds' framework as the denominator identity for the monster Lie algebra. In Eq. \ref{Eq:KNZ},  $q=e^{2\pi i \tau}$, $p=e^{2\pi i \sigma}$ and the $c(n)$'s are the coefficients of the $q$-expansion of $J$.

More details and the history of the moonshine, its proof, applications and developments can be found in \cite{gannon2006monstrous} and \cite{moonshine}. The notes of Borcherds' talk at the 1998 International Congress of Mathematicians have been compiled in \cite{whatismoonshine}, moreover Borcherds has written a short review of the monstrous moonshine in \cite{whatis}.

\section*{Moonshine nowadays}

It turns out that $20$ sporadic groups are \textit{subquotients} of $\M$, and the monstrous moonshine induces moonshines for this \textit{happy family}. This led to a set of conjectures called the \textit{generalised moonshine} by Norton \cite{genmoonnorton}, The proof of which has been completed by Carnahan in 2012 \cite{carnahan2012generalized}.
Robert Griess dubbed the remaining six sporadics which are not subquotients of $\M$, the \textit{pariah groups} \cite{griess}.  Whether these also admit some kind of moonshine was already asked in \cite{conwaynorton}. The pariah O'Nan $O'N$ has been shown to satisfy moonshine properties which connect it to arithmetic \cite{onan} and extend to his subgroup $J_1$ (Janko $1$), another of the pariahs.

Other types of moonshines have also been discovered and studied, such as the Mathieu moonshine \cite{mathieu} which is a correspondence between the elliptic genus of K3 surfaces and the Mathieu sporadic group $M_{24}$. The \textit{umbral moonshine} correspondence \cite{umbralconj} generalizes the Mathieu moonshine and was proved in 2015 \cite{umbralproof}. The reviews \cite{kachru2016elementary} and \cite{anagiannis2018tasi} contain pedagogical introductions to these moonshines.

There are still a lot of unexplained exciting coincidences (also called \textit{moonshines}) which relate sporadic groups to other fields in mathematics. For example, another famous observation of McKay relates the possible orders of the result of the multiplication of two elements of the conjugacy class $2A$ in $\M$, and the Dynkin diagram of the affine algebra $\hat E_8$, as explained in \cite{conway1985simple}. Other hints let one hope for a link between the three largest sporadic groups ($\M$, $\mathbb{B}$ and $\mathrm{Fi}_{24}'$) and the exceptional finite-dimensional semi-simple Lie algebras ($E_8$,$E_7$ and $E_6$) \cite{sporadicandex}. 

Exceptional objects are fascinating, if anything because of their mere existence. Moonshine correspondences hint to deep structures and patterns in mathematics and have already been generous sources of knowledge, though it would be peremptory to say they are fully understood. It is very likely that we are not done being surprised.

\vspace{0.5cm}

I'm thankful to A. Thomas and S. Tornier for their feedback and to R. Duque for his diligent proofreading.

\bibliographystyle{alpha}
\bibliography{ref}

\end{document}